\documentclass[12pt,leqno]{article}
\usepackage{amsmath,amssymb}
\begin{document}
\newtheorem{theorem}{Theorem}[section]
\newtheorem{prop}{Proposition}[section]
\newtheorem{lemma}{Lemma}[section]
\newtheorem{defin}{Definition}[section]
\newtheorem{rem}{Remark}[section]
\newtheorem{example}{Example}[section]
\newtheorem{corol}{Corollary}[section]
\title{Lagrange geometry on tangent manifolds}
\author{{\small by}\vspace{2mm}\\Izu Vaisman}
\date{}
\maketitle
{\def\thefootnote{*}\footnotetext[1]%
{{\it Mathematics Subject Classification:} 53C15, 53C60.
Locally Lagrange metric.}}
\begin{center} \begin{minipage}{12cm}
A{\footnotesize BSTRACT. Lagrange geometry is the geometry of the
tensor field defined by the fiberwise  Hessian of a non degenerate
Lagrangian function on the total space of a tangent bundle.
Finsler geometry is the geometrically most interesting case of
Lagrange geometry. In this paper, we study a generalization, which
consists of replacing the tangent bundle by a general tangent
manifold, and the Lagrangian by a family of compatible, local,
Lagrangian functions. We give several examples, and find the
cohomological obstructions to globalization. Then, we extend the
connections used in Finsler and Lagrange geometry, while giving an
index free presentation of these connections.}
\end{minipage}
\end{center} \vspace{5mm}
\section{Preliminaries} Lagrange geometry
\cite{{Kern},{M1},{M2}} is the extension of Finsler geometry
(e.g., \cite{Bao}) to transversal ``metrics" (non degenerate
quadratic forms) of the vertical foliation (the foliation by
fibers) of a tangent bundle, which are defined as the Hessian of a
non degenerate Lagrangian function. In the present paper, we study
the generalization of Lagrange geometry to arbitrary tangent
manifolds \cite{BC}. The locally Lagrange-symplectic manifolds
\cite{V1} are an important particular case. In this section, we
recall various facts about the geometric structures that we need
for the generalization. Our framework is the $C^\infty$-category,
and we will use the Einstein summation convention, where
convenient.

First, a {\it leafwise locally affine} foliation is a foliation
such that the leaves have a given locally affine structure that
varies smoothly with the leaf. In a different formulation
\cite{V2}, if $M$ is a manifold of dimension $m=p+q$, a
$p$-dimensional locally leafwise affine foliation $ \mathcal{F}$
on $M$ is defined by a maximal, differential, {\it affine atlas}
$\{U_\alpha\}$, with local coordinates $(x^a_\alpha,y^u_\alpha)$
$(a=1,...,q;\,u=1,...,p)$, and transition functions of the local
form
\begin{equation} \label{trfunct} x^a_\beta=x^a_\beta(x^b_\alpha),
\;y^u_\beta=\sum_{v=1}^pA^u_{(\alpha\beta)v}(x^b_\alpha)y^v_\alpha
+B^u_{(\alpha\beta)}(x^b_\alpha)\end{equation} on $U_\alpha\cap
U_\beta$. Then, the leaves of $ \mathcal{F}$ are locally defined
by $x^a=const.$, and their local parallelization is defined by the
vector fields $\partial/\partial y^u$. Furthermore, if the atlas
that defines a locally leafwise affine foliation has a subatlas
such that $B_{(\alpha\beta)}^u=0$ for its transition functions,
the foliation, with the structure defined by the subatlas, will be
called a {\it vector bundle-type foliation}. Notice that, if one
such subatlas exists, similar ones are obtained by coordinate
changes of the local form \begin{equation} \label{translations}
\tilde x^a_\alpha=\tilde x^a_\alpha(x^b_\alpha),\;\tilde y_\alpha
= y^a_\alpha+\xi_{(\alpha\beta)}^a(x^b_\alpha). \end{equation}

For any foliation $ \mathcal{F}$, geometric objects of $M$ that
either project to the space of leaves or, locally, are pull-backs
of objects on the latter are said to be {\it projectable} or {\it
foliated} \cite{{Mol},{V3}}. In particular, a foliated bundle is a
bundle over $M$ with a locally trivializing atlas with foliated
transition functions. The transversal bundle
$\nu\mathcal{F}=TM/T\mathcal{F}$ is foliated. Formulas
(\ref{trfunct}) show that for a locally leafwise affine foliation
$ \mathcal{F}$ the tangent bundles $T\mathcal{F}$, $TM$ are
foliated bundles as well. For a foliated bundle, we can define
foliated cross sections. Notice that, if $ \mathcal{F}$ is a
locally leafwise affine foliation, a vector field on $M$ which is
tangent to $ \mathcal{F}$ is foliated as a vector field, since it
projects to $0$, but, it may not be a foliated cross sections of
$T\mathcal{F}$!

Furthermore, for a locally leafwise affine foliation one also has
{\it leafwise affine} objects, which have an affine character with
respect to the locally affine structure of the leaves. For
instance, a locally leafwise affine function is a function $f\in
C^\infty(M)$ such that $Yf$ is foliated for any local parallel
vector field $Y$ along the leaves of $ \mathcal{F}$. With respect
to the affine atlas, a locally leafwise affine function has the
local expression
\begin{equation} \label{affunct} f=\sum_{u=1}^p
\alpha_u(x^a)y^u+\beta(x^a). \end{equation} A locally leafwise
affine $k$-form is a $k$-form $\lambda$ such that $i(Z)\lambda=0$
for all the tangent vector fields $Z$ of $ \mathcal{F}$ and
$L_Y\lambda$ is a foliated $k$-form for all the parallel fields
$Y$. Then, $\lambda$ has an expression of the form (\ref{affunct})
where $\alpha_u,\beta$ are foliated $k$-forms. A locally
leafwise affine vector field is an infinitesimal automorphism of
the foliation and of the leafwise affine structure, and has the
local expression
\begin{equation}\label{afvect}
X=\sum_{a=1}^q\xi^a(x^b)\frac{\partial}{\partial
x^a}+\sum_{u=1}^p[\sum_{v=1}^p
\lambda^u_v(x^b)y^v+\mu^u(x^b)]\frac{\partial}{\partial y^u}.
\end{equation} Etc. \cite{V2}

Any foliated vector bundle $V\rightarrow M$ produces a sheaf
$\underline{V}$ of germs of differentiable cross sections, and a
sheaf $\underline{V}_{pr}$ of germs of foliated cross sections.
The corresponding cohomology spaces $H^k(M,\underline{V}_{pr})$
may be computed by a de Rham type theorem \cite{V3}. Namely, let
$N\mathcal{F}$ be a complementary (normal) distribution of
$T\mathcal{F}$ in $TM$. The decomposition $TM=N\mathcal{F}\oplus
T\mathcal{F}$ yields a bigrading of differential forms and tensor
fields, and a decomposition of the exterior differential as
\begin{equation} \label{decompd}
d=d'_{(1,0)}+d''_{(0,1)}+\partial_{(2,-1)}.\end{equation} The
operator $d''$ is the exterior differential along the leaves of $
\mathcal{F}$, it has square zero and satisfies the Poincar\'e
lemma. Accordingly, \begin{equation}\label{resol} 0\rightarrow
\underline{V}_{pr}\stackrel{\subseteq}{\rightarrow}\underline{V}_{pr}\otimes_\Phi
\underline{\Omega}^{(0,0)}
\stackrel{d''}{\rightarrow}\underline{V}_{pr}\otimes_\Phi
\underline{\Omega}^{(0,1)}\stackrel{d''}{\rightarrow}...,\end{equation}
where $\Omega$ denotes spaces of differential forms,
$\underline\Omega$ is the corresponding sheaf of differentiable
germs and $\Phi$ is the sheaf of germs of foliated functions,  is
a fine resolution of $ \underline{V}_{pr}$.

Furthermore, if $ \mathcal{F}$ is locally leafwise affine, one
also has the spaces $A^k(M,\mathcal{F})$ of locally leafwise
affine $k$-forms and the corresponding sheaves of germs
$\underline{A} ^k(M,\mathcal{F})$. These sheaves define
interesting cohomology spaces, which may be studied by means of
the exact sequences \cite{V2}
\begin{equation} \label{exact1} 0\rightarrow
\underline{\Omega}^{(k,0)}_{pr}\stackrel{\subseteq}{\rightarrow}
\underline{A}^k(M,\mathcal{F})\stackrel{\pi}{\rightarrow}
\underline{\Omega}^{(k,0)}_{pr}\otimes_\Phi
\underline{T^*\mathcal{F}}_{pr} \rightarrow0,\end{equation} where,
for $f$ defined by (\ref{affunct}),
$\pi(f)=\alpha_u\otimes[dy^u]$, $[dy^u]$ being the projections of
$dy^u$ on $T^*\mathcal{F}$.

It is important to recognize the vector bundle-type foliations
among the locally leafwise affine foliations. First, notice that a
vector bundle-type foliation possesses a global vector field,
which may be seen as the leafwise infinitesimal homothety namely,
\begin{equation} \label{homothety} E=\sum_{u=1}^p
y^u\frac{\partial}{\partial y^u},
\end{equation} called the {\it Euler vector field}. In the general
locally leafwise affine case, (\ref{homothety}) only defines local
vector fields $E_\alpha$ on each coordinate neighborhood
$U_\alpha$, and the differences $E_\beta-E_\alpha$ yield a cocycle
and a cohomology class $[E](\mathcal{F})\in
H^1(M,\underline{T\mathcal{F}}_{pr})$, called the {\it linearity
obstruction} \cite{V2}. It follows easily that the locally
leafwise affine foliation $ \mathcal{F}$ has a vector bundle-type
structure iff $[E](\mathcal{F})=0$ \cite{V2}. With a normal
distribution $N\mathcal{F}$, we may use the foliated version of de
Rham's theorem, and $[E](\mathcal{F})$ will be represented by the
global $T\mathcal{F}$-valued $1$-form $\{d''E_\alpha\}$.
Accordingly, $[E](\mathcal{F})=0$ iff there exists a global vector
field $E$ on $M$, which is tangent to the leaves of $ \mathcal{F}$
and such that $\forall\alpha$,
\begin{equation} \label{genEuler}
E|_{U_\alpha}=E_\alpha+Q_\alpha, \end{equation} where $Q_\alpha$
are projectable. $E$ is defined up to the addition of a global,
projectable, cross section of $T\mathcal{F}$, and these vector
fields $E$ will be called {\it Euler vector fields}. The choice of
a Euler vector field $E$ is equivalent with the choice of the
vector bundle-type structure of the foliation.

We also recall the following result \cite{V2}: the vector
bundle-type foliation $ \mathcal{F}$ on $M$ is a vector bundle
fibration $M\rightarrow N$ iff the leaves are simply connected and
the flat connections defined by the locally affine structure of
the leaves are complete. \begin{example}\label{example0} {\rm On
the torus $ \mathbb{T}^{p+q}$ with the Euclidean coordinates
$(x^a,y^u)$ defined up to translations $$\tilde
x^a=x^a+h^a,\,\tilde y^u=y^u+k^u,\hspace{5mm}h^a,k^u\in
\mathbb{Z},$$ the foliation $x^a=const.$ is locally leafwise
affine and has the normal bundle $dy^u=0$. The linearity
obstruction $[E]$ is represented by the form $\sum_{u=1}^q
dy^u\otimes(\partial/\partial y^u)$, which is not $d''$-exact.
Therefore, $[E]\neq0$, and $\mathcal{F}$ is not a vector
bundle-type foliation.}\end{example}
\begin{example}\label{example1} {\rm Consider the compact nilmanifold
$M(1,p)=\Gamma(1,p)\backslash H(1,p)$ where
\begin{equation} \label{Heisenberg}
H(1,p)=\left\{\left(\begin{array}{ccc} Id_p&X&Z\\
0&1&y\\0&0&1 \end{array}\right)\,/\,X,Z\in \mathbb{R}^p,\,y\in
\mathbb{R}\right\}\end{equation} is the generalized Heisenberg
group, and $\Gamma(1,p)$ is the subgroup of matrices with integer
entries. $M(1,p)$ has an affine atlas with the transition
functions
\begin{equation} \label{trHeis} \tilde x^i=x^i+a^i,\,\tilde
y=y+b,\, \tilde z^i=z^i+a^iy+c^i,\end{equation} where $x^i,z^i$
$(i=1,...,p)$ are the entries of $X,Z$, respectively, and
$a^i,b,c^i$ are integers.  Accordingly, the local equations
$x^i=const.,\, y=const.$ define a  locally leafwise affine
foliation $ \mathcal{F}$ of $M$, which, in fact, is a fibration by
$p$-dimensional tori over a $(p+1)$-dimensional torus. The
manifold $M$ is parallelizable by the global vector fields
\begin{equation}\label{paralHeis}
\frac{\partial}{\partial x^i}, \frac{\partial}{\partial
y}+\sum_{i=1}^px^i\frac{\partial}{\partial
z^i},\frac{\partial}{\partial z^i},\end{equation} and the global
$1$-forms
\begin{equation}\label{paralHeis2} dx^i,dz^i-x^idy,dy,
\end{equation} and we see that $$span\left\{\frac{\partial}{\partial x^i},
\frac{\partial}{\partial
y}+\sum_{i=1}^px^i\frac{\partial}{\partial z^i}\right\}$$ may
serve as a normal bundle of $ \mathcal{F}$. It follows that the
linearity obstruction is represented by
$$\sum_{i=1}^p (dz^i-x^idy)\otimes\frac{\partial}{\partial z^i},$$
which is not $d''$-exact. Therefore, $ \mathcal{F}$ is not a
vector bundle-type foliation.}
\end{example}
\begin{example}\label{example2} {\rm
Take the {\it real Hopf manifold} $H^{(p+q)}=S^{p+q-1}\times
S^1$ seen as $ (\mathbb{R}^q\times
\mathbb{R}^p\,\backslash\,\{0\})\,/\,G_\lambda$, where
$\lambda\in(0,1)$ is constant and $G_\lambda$ is the group
\begin{equation}\label{groupG} \tilde x^a=\lambda^nx^a,
\tilde y^u=\lambda^ny^u,\hspace{5mm}n\in\mathbb{Z},\end{equation}
where $x^a,y^u$ are the natural coordinates of $ \mathbb{R}^q$ and
$ \mathbb{R}^p$, respectively. Then, the local equations
$x^a=const.$ define a vector bundle-type foliation, which has the
global Euler field $E=\sum_{u=1}^qy^u(\partial/\partial y^u)$.
This example shows that compact manifolds may
have vector bundle-type foliations.}\end{example}
\begin{example}\label{example3} {\rm
Consider the manifold
\begin{equation}\label{example31}
M^{2n}=[( \mathbb{R}^n\,\backslash\,\{0\})\times \mathbb{R}^n]
\,/\,K_\lambda,\end{equation} where $\lambda\in(0,1)$ and
$K_\lambda$ is the cyclic group generated by the transformation
\begin{equation} \label{eqexample3} \tilde x^i=\lambda x^i,\,
\tilde y^i=\lambda
y^i+(1-\lambda)\frac{x^i}{\sqrt{\sum_{j=1}^n(x^j)^2}}
\end{equation} $(i=1,...n)$. It is easy to check that the equality
$$E=\sum_{i=1}^n y^i
\frac{\partial}{\partial
y^i}-\sum_{i=1}^n\frac{x^i}{\sqrt{\sum_{j=1}^n(x^j)^2}}\frac{\partial}{\partial
y^i}$$ defines a global vector field on $M$, which has the
property of the Euler field for the foliation $x^i=const.$
Therefore, the latter is a vector bundle-type foliation. The
change of coordinates
$$x'^i=x^i,\,y'^i=y^i-\frac{x^i}{\sqrt{\sum_{j=1}^n(x^j)^2}}$$
provides a vector bundle-type atlas, and (\ref{eqexample3})
becomes
$$\tilde x'^i=\lambda x^i,\,\tilde y'^i=\lambda y^i.$$
This shows that $M$ is the tangent bundle of the Hopf manifold
$H^n$ defined in Example \ref{example2}}.
\end{example}

Now, let us recall the basics of tangent manifolds \cite{BC}. An
{\it almost tangent structure} on a manifold $M$ is a tensor field
$S\in \Gamma End(TM)$ such that \begin{equation}\label{almosttg}
S^2=0,\;im\,S=ker\,S.\end{equation} In particular, the dimension
of $M$ must be even, say $2n$, and $rank\,S=n$. Furthermore, $S$
is a {\it tangent structure} if it is integrable i.e., locally,
$S$ looks like the vertical twisting homomorphism of a tangent
bundle. This means that there exists an atlas with local
coordinate $(x^i,y^i)$ $(i=1,...,n)$ such that
\begin{equation}\label{strtg} S\left(\frac {\partial}{\partial
x^i}\right)=\frac{\partial}{\partial y^i},\, S\left(\frac
{\partial}{\partial y^i}\right)=0. \end{equation} The
integrability property is equivalent with the annulation of the
Nijenhuis tensor \begin{equation} \label{NijS} \mathcal{N}_S(X,Y)=
[SX,SY]-S[SX,Y]-S[X,SY]+S^2[X,Y]=0.
\end{equation}
A pair $(M,S)$, where $S$ is a tangent structure, is
called a {\it tangent manifold}.

On a tangent manifold $(M,S)$, the distribution $im\,S$ is
integrable, and defines the {\it vertical foliation} $
\mathcal{V}$ with $T\mathcal{V}=im\,S$. It is easy to see that the
transition functions of the local coordinates of (\ref{strtg}) are
of the local form (\ref{trfunct}) with $q=p=n$ and
\begin{equation}\label{transtg}
A^i_{(\alpha\beta)j}=\frac{\partial x^i_\beta}{\partial
x^j_\alpha}. \end{equation} Therefore, $ \mathcal{V}$ is a locally
leafwise affine foliation, and the local parallel vector fields
along the leaves are the vector fields of the form $SX$, where $X$
is a foliated vector field. In particular, a tangent manifold has
local Euler fields $E_\alpha$, and a linearity obstruction $[E]\in
H^1(M, \underline{T\mathcal{V}}_{pr})$. If $[E]=0$, the foliation
$ \mathcal{V}$ will be a vector bundle-type foliation, and $M$ has
global Euler vector fields $E$ defined up to the addition of a
foliated cross section of $T\mathcal{V}$. Furthermore, if we fix
the vector-bundle type structure by fixing a Euler vector field
$E$, the triple $(M,S,E)$ will be called a {\it bundle-type
tangent manifold}.

Using the general result of \cite{V2}, we see that a tangent
manifold is a tangent bundle iff it is a bundle-type tangent
manifold and the vertical foliation has simply connected, affinely
complete leaves.
\begin{example}\label{example4}{\rm The Hopf manifold $H^{2n}$ of
Example \ref{example2} with $q=p=n$ and $S$ defined by
(\ref{strtg}) is a compact, bundle-type, tangent
manifold.}\end{example}
\begin{example}\label{example5} {\rm The torus of Example
\ref{example0} with $q=p$ and $S$ of (\ref{strtg}) is a compact
non bundle-type tangent manifold.}\end{example}
\begin{example}\label{example6} {\rm The manifold $M(1,p)\times(
\mathbb{R}/\mathbb{Z})$, with the coordinates of Example
\ref{example1} and a new coordinate $t$ on $ \mathbb{R}$, and with
$S$ defined by
\begin{equation} \label{ex1tg} S\left(\frac
{\partial}{\partial x^i}\right)=\frac{\partial}{\partial z^i},\,
S\left(\frac {\partial}{\partial
y}\right)=\frac{\partial}{\partial t},\, S\left(\frac
{\partial}{\partial z^i}\right)=0,\, S\left(\frac
{\partial}{\partial t}\right)=0
\end{equation}  is a compact non bundle-type tangent
manifold. The linearity obstruction $[E]$ of this manifold is
represented by
$$\sum_{i=1}^p (dz^i-x^idy)\otimes\frac{\partial}{\partial z^i}
+dt\otimes\frac{\partial}{\partial t},$$ and $[E]\neq0$.}
\end{example}

Tangent bundles posses {\it second order vector fields} ({\it
semisprays} in \cite{M1}), so called because they may be locally
expressed by a system of second order, ordinary, differential
equations. A priori, such vector fields may be defined on any
tangent manifold \cite{V4} namely, the vector field $X\in \Gamma
TM$ ($\Gamma$ denotes the space of global cross sections) is of
the {\it second order} if $SX|_{U_\alpha}-E_\alpha$ is foliated
for all $\alpha$. But, this condition means that $SX$ is a global Euler
vector field, hence, only the bundle-type tangent manifolds can
have global second order vector fields.

It is important to point out that, just like on tangent bundles
(e.g., \cite{{Kern},{M1},{V5}}), if $(M,S,E)$ is a bundle-type
tangent manifold, and $X$ is a second order vector field on $M$,
the Lie derivative $F=L_XS$ defines an almost product structure on
$M$ ($F^2=Id$), with the associated projectors
\begin{equation}
\label{almprod}V=\frac{1}{2}(Id+F),\;H=\frac{1}{2}(Id-F),\end{equation}
such that $im\,V=T\mathcal{V}$ and $im\,H$ is a normal
distribution $N\mathcal{V}$ of the vertical foliation $
\mathcal{V}$.

Finally, we give \begin{defin}\label{tangautomorf} {\rm A vector field
$X$ on a tangent manifold $(M,S)$ is a {\it tangential
infinitesimal automorphism} if $L_XS=0$ ($L$ denotes the Lie
derivative).} \end{defin}

Obviously, a tangential infinitesimal automorphism $X$ preserves
the foliation $ \mathcal{V}$ and its leafwise affine structure.
Therefore, $X$ is a leafwise affine vector field with respect to $
\mathcal{V}$. Furthermore, in the bundle-type case, if $E$ is a
Euler vector field, $[X,E]$ is a foliated cross section of
$T\mathcal{V}$.
\section{Locally Lagrange spaces} Lagrange geometry is motivated
by physics and, essentially, it is the study of geometric objects
and constructions that are transversal to the vertical foliation
of a tangent bundle and are associated with a {\it Lagrangian} (a
name taken from Lagrangian mechanics) i.e., a function on the
total space of the tangent bundle. (See \cite{M1} and the
$d$-objects defined there.) Here, we use the same approach for a
general tangent manifold $(M,S)$, and we refer to functions on $M$
as {\it global Lagrangians} and to functions on open subsets as
{\it local Lagrangians}.

If $ \mathcal{L}$ is a Lagrangian, the derivatives in the vertical
directions yield symmetric tensor fields of $M$ defined by
\begin{equation} \label{Hessk} (Hess_{(k)}\mathcal{L})_x(X_1,...,X_k)
=(S\tilde X_k)\cdots(S\tilde X_1)\mathcal{L}|_x,\;x\in M,\; X_i\in
T_xM,\end{equation} where $\tilde X_i$ $(i=1,...,k)$ are
extensions of $X_i$ to local, $ \mathcal{V}$-foliated, vector
fields on $M$. (Of course, the result does not depend on the
choice of the extensions $\tilde X_i$.) $Hess_k\mathcal{L}$ is
called the $k$-{\it  Hessian} of $ \mathcal{L}$. Notice that
definition (\ref{Hessk}) may also be replaced by the recurrence
formula \begin{equation}\label{Hessrec}
(Hess_{(k)}\mathcal{L})_x(\tilde X_1,...,\tilde X_k)= [L_{S\tilde
X_k}(Hess_{k-1} \mathcal{L})]_x(\tilde X_1,...,\tilde X_{k-1}),
\end{equation} where the arguments are foliated vector fields.

It is worthwhile to notice the following general property
\begin{prop} \label{HessX} for any function $ \mathcal{L}\in
C^\infty(M)$, any tangential infinitesimal automorphism $X$ of the
tangent manifold $(M,S)$, and any $k=1,2,...$, one has
\begin{equation} \label{LieHessX}
Hess_k(X\mathcal{L})=L_X(Hess_k\mathcal{L}). \end{equation}
\end{prop} \noindent{\bf Proof.} Proceed by induction on $k$, while
evaluating the Hessian of $X\mathcal{L}$ on foliated arguments and
using the recurrence formula (\ref{Hessrec}). Q.e.d.

For $k=1$, we get a $1$-form, say $\theta_\mathcal{L}$, and for
$k=2$, we get the usual Hessian of $ \mathcal{L}$ with respect to
the affine vertical coordinates $y^i$ (see Section 1), hereafter
to be denoted by either $Hess\,\mathcal{L}$ or $g_\mathcal{L}$.
Obviously, $g_\mathcal{L}$ vanishes whenever one of the arguments
is vertical hence, it yields a well defined cross section of the
symmetric tensor product $\odot^2\nu^*\mathcal{V}$
$(\nu\mathcal{V}=TM/T\mathcal{V})$, which we continue to denote by
$g_\mathcal{L}$. If $g_\mathcal{L}$ is non degenerate on the
transversal bundle $\nu\mathcal{F}$, the Lagrangian $ \mathcal{L}$
is said to be {\it regular} and $g_\mathcal{L}$ is called a {\it
(local) Lagrangian metric}. We note that if the domain of $
\mathcal{L}$ is connected, the regularity of $ \mathcal{L}$ also
implies that $g_\mathcal{L}$ is of a constant signature. With
respect to the local coordinates of (\ref{strtg}), one has
\begin{equation}\label{Hesslocal} \theta_{ \mathcal{L}}=
\frac{\partial\mathcal{L}}{\partial y^i}dx^i,\;g_{ \mathcal{L}}=
\frac{1}{2}\frac{\partial^2\mathcal{L}}{\partial y^i\partial
y^j}dx^i\odot dx^j.\end{equation}

Lagrangian mechanics shows the interest of one more geometric
object related to a Lagrangian namely, the differential $2$-form
\begin{equation}\label{formasimpl}
\omega_\mathcal{L}=d\theta_\mathcal{L}=
\frac{\partial^2\mathcal{L}}{\partial x^i\partial y^j}dx^i \wedge
dx^j + \frac{\partial^2\mathcal{L}}{\partial y^i\partial
y^j}dy^i\wedge dx^j.\end{equation} If $ \mathcal{L}$ is a regular
Lagrangian $\omega_\mathcal{L}$ is a symplectic form, called the
{\it Lagrangian symplectic form}.

In \cite{{V1},{V4}}, we studied particular symplectic forms
$\Omega$ on a tangent manifold $(M,S)$ that are {\it compatible
with the tangent structure} $S$ in the sense that
\begin{equation} \label{symplcomp}
\Omega(X,SY)=\Omega(Y,SX).\end{equation} If this happens, $\Omega$
is called a {\it locally Lagrangian-symplectic form} since the
compatibility property is equivalent with the existence of an open
covering $M=\cup U_\alpha$, and of local regular Lagrangian
functions $ \mathcal{L}_\alpha$ on $U_\alpha$, such that,
$\Omega|_{U_\alpha}=\omega_{ \mathcal{L}_\alpha}$ for all
$\alpha$. On the intersections $U_\alpha\cap U_\beta$ the local
Lagrangians satisfy a compatibility relation of the form
\begin{equation} \label{Lrel}
\mathcal{L}_\beta-\mathcal{L}_\alpha=a(\varphi_{(\alpha\beta)})
+b_{(\alpha\beta)},\end{equation} where $\varphi_{(\alpha\beta)}$
is a closed, foliated $1$-form, $b_{(\alpha\beta)}$ is a foliated
function, and $a(\varphi)=\varphi_iy^i$ where the local
coordinates and components are taken either in $U_\alpha$ or in
$U_\beta$. Furthermore, if it is possible to find a compatible (in
the sense of (\ref{Lrel})) global Lagrangian $ \mathcal{L}$,
$\Omega$ is a {\it global Lagrangian symplectic form}. Conditions
for the existence of a global Lagrangian were given in
\cite{{V1},{V4}}. In particular, a globally Lagrangian-symplectic
manifold $M^{2n}$ cannot be compact since it has the exact volume
form $\omega_{ \mathcal{L}}^n$.

Following the same idea, we give \begin{defin} \label{varlocLagr}
{\rm Let $(M^{2n},S)$ be a tangent manifold, and
$g\in\Gamma\odot^2\nu^*\mathcal{V}$ a non degenerate tensor field.
Then $g$ is a {\it locally Lagrangian metric (structure)} on $M$
if there exists an open covering $M=\cup U_\alpha$ with local
regular Lagrangian functions $ \mathcal{L}_\alpha$ on $U_\alpha$
such that $g|_{U_\alpha}=g_{ \mathcal{L}_\alpha}
=Hess\,\mathcal{L}_\alpha$ for all $\alpha$. The triple $(M,S,g)$
will be called a {\it locally Lagrange space or
manifold}.}\end{defin}

It is easy to see that the local Lagrangians $ \mathcal{L}_\alpha$
of a locally Lagrange space must again satisfy the compatibility
relations (\ref{Lrel}), where the $1$-forms
$\varphi_{(\alpha\beta)}$ may not be closed. In particular, we see
that a locally Lagrangian-symplectic manifold is a locally
Lagrange space with the metric defined by \cite{V1}
\begin{equation}\label{metricoflocally Lagrangian-symplectic}
g([X],[Y])=\Omega(SX,Y),\end{equation} where $X,Y\in\Gamma TM$ and
$[X],[Y]$ are the corresponding projections on $\nu\mathcal{F}$.
Furthermore, if there exists a global Lagrangian $ \mathcal{L}$
that is related by (\ref{Lrel}) with the local Lagrangians of the
structure, $(M,S,g,\mathcal{L})$ will be a {\it globally Lagrange
space}. A globally Lagrange space also is a globally
Lagrangian-symplectic manifold hence, it cannot be compact.

We can give a global characterization of the locally Lagrange
metrics. First, we notice that the bundles
$\otimes^k\nu^*\mathcal{V}$ of covariant tensors transversal to
the vertical foliation $ \mathcal{V}$ of a tangent manifold
$(M,S)$ may also be seen as the bundles of covariant tensors on
$M$ that vanish if evaluated on arguments one of which belongs to
$im\,S$. (This holds because $\nu*\mathcal{V}\subseteq T^*M$.) In
particular, a transversal metric $g$ of $ \mathcal{V}$ may be seen
as a symmetric $2$-covariant tensor field $g$ on $M$ which is
annihilated by $im\,S$. With $g$, one associates a $3$-covariant
tensor, called the {\it derivative} or {\it Cartan tensor}
\cite{{Bao},{M1},{M2}} defined by
\begin{equation}\label{defC} C_x(X,Y,Z)=(L_{S\tilde X}g)_x
(Y,Z),\;x\in M,\;X,Y,Z\in T_xM,\end{equation} where $\tilde X$ is
a foliated extension of $X$. Obviously,
$C\in\Gamma\otimes^3\nu^*\mathcal{V}$. Then, we get
\begin{prop}\label{propos25} The transversal metric $g$ of the
vertical foliation $ \mathcal{V}$ of a tangent manifold $(M,S)$ is
a locally Lagrange metric iff the tensor field $C$ is totally
symmetric. \end{prop} \noindent {\bf Proof.} Since
$$C_{ijk}=C(\frac{\partial}{\partial x^i},\frac{\partial}{\partial x^j},
\frac{\partial}{\partial x^k})=\frac{\partial g_{jk}}{\partial
y^i},$$ the symmetry of $C$ is equivalent with the existence of
the required local Lagrangians $ \mathcal{L}$. Q.e.d.

We give a number of examples of locally Lagrange manifolds.
\begin{example}\label{example21} {\rm Consider the torus of Example
\ref{example5}. Then $$ \mathcal{L}=\frac{1}{2}\sum_{i=1}^n(y^i)^2
$$ define compatible local Lagrangians with the corresponding
Lagrange metric $\sum_{i=1}^n(dx^i)^2$. (Notice also the existence
of the locally Lagrange symplectic form $\Omega=\sum_{i=1}^{n}dx^i
\wedge dy^i$.)}
\end{example}
\begin{example}\label{example22}{\rm Consider the tangent manifold
$M(1,p)\times( \mathbb{R}/\mathbb{Z})$ of Example \ref{example6},
with the tangent structure defined by (\ref{ex1tg}). The $
\mathcal{V}$-transversal metric
$$\sum_{i=1}^p(dx^i)^2+(dy)^2$$ is the Lagrange metric of the
local compatible Lagrangians
$$\frac{1}{2}(\sum_{i=1}^p(z^i)^2+t^2).$$ (In this example the forms
$\varphi_{(\alpha\beta)}$ of (\ref{Lrel}) are not closed.)}
\end{example}

Examples \ref{example21}, \ref{example22} are interesting because
the manifolds involved are compact manifolds.
\begin{example}\label{example23} {\rm The manifold $M^{2n}$ of
Example \ref{example3} is diffeomorphic with the tangent bundle
$TH^n$. With the coordinates $(x'^i,y'^i)$ (see Example
\ref{example3}), we see that the function $$
\mathcal{L}=\frac{\sum_{i=1}^{n}(y'^i)^2}{2\sum_{i=1}^{n}(x'^i)^2}$$
is a global, regular Lagrangian, and it produces a positive
definite Lagrange metric.}\end{example}
\begin{example}\label{example24}
{\rm Consider the Hopf manifold $H^{2n}$ of Example \ref{example4}
with the tangent structure (\ref{strtg}), and define the local
compatible Lagrangians \begin{equation}\label{metricaln}
\mathcal{L}=\frac{1}{2}\ln{\rho},\;\;\rho=\sum_{i=1}^n[(x^i)^2+(y^i)^2].
\end{equation} An easy computation yields
\begin{equation} \label{Hessln}
\frac{\partial^2\mathcal{L}}{\partial y^i\partial y^j}
=-\frac{1}{\rho^2}(2y^iy^j-\rho\delta_{ij}).\end{equation} The
determinant of the Hessian (\ref{Hessln}) can be easily computed
as a characteristic polynomial and we get
$$det\left(\frac{\partial^2\mathcal{L}}{\partial y^i\partial y^j}\right)=
\frac{\sum_{i=1}^n[(x^i)^2-(y^i)^2]}
{\{\sum_{i=1}^n[(x^i)^2+(y^i)^2]\}^{n+1}}.$$ Now, the local
equation $$\sum_{i=1}^n(x^i)^2=\sum_{i=1}^n(y^i)^2$$ defines a
global hypersurface $\Sigma$ of $H^{2n}$, and (\ref{Hessln})
provides a locally Lagrange metric structure on
$H^{2n}\backslash\Sigma$.}
\end{example}
\begin{example} \label{example26} {\rm On any tangent manifold $(M,S)$,
any non degenerate, foliated, transversal metric $g$ of the
vertical foliation (if such a metric exists \cite{Mol}) is locally
Lagrange. Indeed, this kind of metric is characterized by $C=0$,
and the result follows from Proposition
\ref{propos25}.}\end{example}

A natural question implied by Definition \ref{varlocLagr} is:
assume that $(M,S,g,\mathcal{L}_\alpha)$ is a locally Lagrange
space; what conditions ensure the existence of a global
compatible, regular Lagrangian?

The compatibility relations (\ref{Lrel}) endow $M$ with an
$\underline{A}^0$-valued $1$-cocycle defined by any of the members
of equation (\ref{Lrel}), hence, with a cohomology class $
\mathcal{G}\in H^1(M,\underline{A}^0)$, which we call the {\it
total Lagrangian obstruction}, and it is obvious that $
\mathcal{G}=0$ iff the manifold $M$ with the indicated structure
is a globally Lagrange space.

Furthermore, the total Lagrangian obstruction may be decomposed
into two components determined by the exact sequence
(\ref{exact1}) with $k=0$, which in our case becomes
\begin{equation} \label{exact3}
0\rightarrow\Phi\stackrel{\subseteq}{\rightarrow}\underline{A}^0
(M,\mathcal{V})\stackrel{\pi'}{\rightarrow}\underline{\Omega}^{(1,0)}_{pr}
\rightarrow0,\end{equation} where $\pi'$ is the composition of the
projection $\pi$ of (\ref{exact1}) by $S$.

It is easy to see that the connecting homomorphism of the exact
cohomology sequence of (\ref{exact3}) is zero in dimension $0$.
Accordingly, we get the exact sequence
\begin{equation} \label{exactobstruction}
0\rightarrow
H^1(M,\Phi)\stackrel{\iota*}{\rightarrow}H^1(M,\underline{A}^0
)\stackrel{\pi*}{\rightarrow}H^1(M,\underline{\Omega}^{(1,0)}_{pr})
\stackrel{\partial}{\rightarrow}H^2(M,\Phi)\rightarrow\cdots,\end{equation}
where $\iota^*,\pi^*$ are induced by the inclusion and the
homomorphism $\pi'$ of (\ref{exact3}). Accordingly, we get the
cohomology class $ \mathcal{G}_1=\pi^*( \mathcal{G})\in
H^1(M,\underline{\Omega}^{(1,0)}_{pr})$, and we call it the {\it
first Lagrangian obstruction}. $ \mathcal{G}_1=0$ is a necessary
condition for $M$ to be a globally Lagrange space. Furthermore, if
$ \mathcal{G}_1=0$, the exact sequence (\ref{exactobstruction})
tells us that there exist a unique cohomology class $
\mathcal{G}_2\in H^1(M,\Phi)$ such that $ \mathcal{G}=\iota^*(
\mathcal{G}_2)$. We call $ \mathcal{G}_2$ the {\it second
Lagrangian obstruction} of the given structure, and $
\mathcal{G}=0$ iff $ \mathcal{G}_1=0$ and $ \mathcal{G}_2=0$.

We summarize the previous analysis in \begin{prop}
\label{propos21} The locally Lagrange space
$(M,S,g,\mathcal{L}_\alpha)$ is a globally Lagrange space iff both
the first and the second Lagrangian obstructions exist and are
equal to zero. \end{prop}

Let us assume that a choice of a normal bundle $N\mathcal{V}$ has
been made. Then we can use the de Rham theorem associated with the
relevant resolution (\ref{resol}) in order to get a representation
of the Lagrangian obstructions. The definition of $ \mathcal{G}_1$
shows that the first Lagrangian obstruction is represented by the
cocycle $\{\theta_{ \mathcal{L}_\beta}-\theta_{
\mathcal{L}_\alpha}\}$. Accordingly, $ \mathcal{G}_1$ may be seen
as the $d''$-cohomology class of the global form $\Theta$ of
type $(1,1)$ defined by gluing up the local forms $\{d''\theta_{
\mathcal{L}_\alpha}\}$. If we follow the notation of \cite{V3} and
take bases \begin{equation}  \label{bases}
N\mathcal{V}=span\left\{X_i=\frac{\partial}{\partial x^i}-t^j_i
\frac{\partial}{\partial y^j}\right\},\;T\mathcal{V}=span\left\{
Y_i=\frac{\partial}{\partial y^i}\right\},\end{equation} with the
dual cobases \begin{equation}\label{cobases} \begin{array}{l}
N^*\mathcal{V}=ann(T\mathcal{V})=span\{dx^i\},\vspace{1mm}\\
T^*\mathcal{V}=ann(N
\mathcal{V})=span\{\vartheta^i=dy^i+t^i_jdx^j\},\end{array}\end{equation}
where $t^i_j(x^i,y^i)$ are local functions, we get
\begin{equation}\label{formagamma}\Theta=
\frac{\partial^2 \mathcal{L}_\alpha}{\partial y^i\partial
y^j}\vartheta^i\wedge dx^j.\end{equation} The result may be
written as
\begin{prop} \label{propos22} Let $(M,S,g,\mathcal{L}_\alpha)$
be a locally Lagrange space. Then, each choice of a normal bundle
$N \mathcal{V}$ defines an almost symplectic structure of $M$,
given by the non degenerate $d''$-closed $2$-form $\Theta$. The
first Lagrangian obstruction $ \mathcal{G}_1$ vanishes iff the
form $\Theta$ is $d''$-exact.
\end{prop} \begin{corol}\label{corolar20} A compact, connected,
bundle-type, tangent manifold, with the Euler vector field $E$ has
no locally Lagrange metric $g$ such that $L_Eg=sg$ where $s$ is a
function that never takes the value $-1$.
\end{corol} \noindent{\bf Proof.} Essentially, the hypothesis on
$E$ means $E$ cannot be a conformal infinitesimal automorphism of
$g$. From (\ref{formagamma}) we get
\begin{equation}\label{eqPsi} \Psi=\frac{1}{n!}\Theta^n=
(-1)^{\frac{n(n+1)}{2}}det\left(\frac{\partial^2
\mathcal{L}_\alpha}{\partial y^i\partial y^j}\right)
\end{equation} $$\cdot dx^1\wedge\ldots\wedge dx^n\wedge dy^1\wedge\ldots
\wedge dy^n$$ and
\begin{equation} \label{Liegamman} L_E\Psi=(-1)^{\frac{n(n+1)}{2}}\left[E\,
det\left(\frac{\partial^2 \mathcal{L}_\alpha}{\partial y^i\partial
y^j}\right)+n\,det\left(\frac{\partial^2
\mathcal{L}_\alpha}{\partial y^i\partial y^j}\right)\right]
\end{equation} $$\cdot dx^1\wedge\ldots\wedge dx^n\wedge dy^1\wedge\ldots
\wedge dy^n,$$ where the local coordinates belong to an affine
atlas where $E=y^i(\partial/\partial y^i)$. If $M$ is compact,
$\int_M L_E\Psi=0$, and the coefficient of the right hand side of
(\ref{Liegamman}) cannot have a fixed sign. But, the latter property holds
under the hypothesis of the corollary. Q.e.d.

For instance, the Hopf manifold $H^n$ has no locally Lagrange
metric with homogeneous with respect to the coordinates $(y^i)$
Lagrangians $ \mathcal{L}_\alpha$. Indeed, homogeneity of degree
$s\neq-1$ is impossible because of the previous corollary, and
homogeneity of degree $-1$ contradicts the transition relations
(\ref{Lrel}).
\begin{rem}\label{remexHopf} {\rm Because of Corollary
\ref{corolar20}, we conjecture that a compact, bundle-type,
tangent manifold cannot have a locally Lagrange metric.}
\end{rem}
\begin{prop}\label{corolar21} The first Lagrangian obstruction
of a locally Lagrange metric structure of $M$ with the local
Lagrangians $\{ \mathcal{L}_\alpha\}$ vanishes iff there exists a
subordinated structure $\{\tilde{ \mathcal{L}}_\alpha\}$ such that
the $1$-forms $\theta_{\tilde{\mathcal{L}}_\alpha}$ glue up to a
global $1$-form. This subordinated structure defines a locally
Lagrangian-symplectic structure on the manifold $M$. Furthermore,
in this case the second Lagrangian obstruction $ \mathcal{G}_2$ is
represented by the global $d''$-closed form $\kappa$ of type
$(0,1)$ defined by gluing up the local forms $\{ d''\tilde{
\mathcal{L}}_\alpha\}$.\end{prop} \noindent{\bf Proof.} Under the
hypothesis, there exists a global form $\lambda$ of type $(1,0)$
such that $\Theta=d''\theta_{ \mathcal{L}_\alpha}=d''\lambda$,
therefore, $\theta_{
\mathcal{L}_\alpha}=\lambda|_{U_\alpha}+\xi_\alpha$, with some
local foliated $1$-forms $\xi_\alpha=\xi_{\alpha,i}(x^j)dx^i$.
Accordingly, we get
\begin{equation}\label{eqcorol}
\frac{\partial( \mathcal{L}_\beta-\mathcal{L}_\alpha)}{\partial y^i}
=\xi_{\beta,i}-\xi_{\alpha,i}, \end{equation} whence $$
\mathcal{L}_\beta-
\mathcal{L}_\alpha=a(\xi_\beta-\xi_\alpha)+b_{(\alpha\beta)},$$
where $a$ has the same meaning as in (\ref{Lrel}) and
$b_{(\alpha\beta)}$ are foliated functions. Now, if we define
\begin{equation} \label{eq2corol}
\tilde{\mathcal{L}}_\alpha=\mathcal{L}_\alpha-a(\xi_\alpha)
\end{equation} we are done. The last assertion follows from the definition
of $ \mathcal{G}_2$. Q.e.d.
\begin{corol} \label{corolar211}
The locally Lagrange metric of Proposition \ref{corolar21} is
defined by a global Lagrangian iff $\kappa=d''k$ for a function
$k\in C^\infty(M)$. \end{corol}

In order to give an application of this result we recall
\begin{lemma}\label{contractibility} For the vertical foliation $
\mathcal{V}$ of a tangent bundle $TN$, one has $H^k(TN,\Phi)=0$
for any $k>0$.\end{lemma} \noindent{\bf Proof.} Use a normal
bundle $N\mathcal{V}$, and let $\lambda$ be a $d''$-closed form of
type $(p,q)$ on $TN$. Since the fibers of $TN$ are contractible,
if $N=\cup U_\alpha$ is a covering by small enough,
$TN$-trivializing neighborhoods, we have
$\lambda|_{p^{-1}(U_\alpha)}=d''\mu_\alpha$ ($p:TN\rightarrow N$)
for some local forms $\mu_\alpha$ of type $(p,q-1)$. The local
forms $\mu_\alpha$ can be glued up to a global form $\mu$ by means
of the pullback to $TN$ of a partition of unity on $N$, i.e.,  by
means of foliated functions. Accordingly, we will have
$\lambda=d''\mu$. Q.e.d.

From Corollary \ref{corolar211} and Lemma \ref{contractibility} we
get
\begin{prop}\label{propos23} Any locally Lagrange metric of a
tangent bundle $TN$ is a globally Lagrange metric.
\end{prop}
\begin{rem} \label{observatie} {\rm
Propositions \ref{propos25}, \ref{propos23} imply that, in the
case of a tangent bundle $M=TN$, the symmetry of $C$ is a
necessary and sufficient condition for $g$ to be a global
Lagrangian metric. It was well known that this condition is
necessary \cite{M1}. On the other hand the metrics of \cite{M1}
usually are differentiable only on the complement of the zero
section of $TN$, where Proposition \ref{propos23} does not hold,
hence, the condition is not a sufficient one.} \end{rem}

We also mention the inclusion
$\sigma:\underline{Z}^{(1,0)}_{pr}\rightarrow\underline{\Omega}^{(1,0)}_{pr}$,
where $Z$ denotes spaces of closed forms, and the obvious
\begin{prop}\label{propos24} The locally Lagrange metric structure
defined by $\{ \mathcal{L}_\alpha\}$ is reducible to a locally
Lagrangian-symplectic structure iff $ \mathcal{G}_1\in
im\,\sigma^*$, where $\sigma^*$ is induced by $\sigma$ in
cohomology. \end{prop}

Other important notions are defined by \begin{defin}
\label{automorfLagr} {\rm Let $(M,S,g)$ be a locally Lagrange
space, and $X\in\Gamma TM$. Then: i) $X$ is a {\it Lagrange
infinitesimal automorphism} if $L_Xg=0$, where $g$ is seen as a
$2$-covariant tensor field on $M$; ii) $X$ is a {\it strong
Lagrange infinitesimal automorphism} if it is a Lagrange and a
tangential infinitesimal automorphism of $(M,S)$, simultaneously.}
\end{defin}

Notice that \begin{equation} \label{eqptderivLie}
(L_Xg)(Y,SZ)=-g(Y,[X,SZ])\hspace{5mm}(X,Y,Z\in\Gamma TM).
\end{equation}
From (\ref{eqptderivLie}) and the non degeneracy of $g$ on $\nu\mathcal{V}$
it follows that a Lagrange infinitesimal automorphism necessarily is a
$\mathcal{V}$-projectable vector field. But, it may not be locally leafwise affine.
Indeed, if $g$ is a foliated metric of $\nu\mathcal{V}$
(Example \ref{example26}) every tangent vector field of $\mathcal{V}$ is a
Lagrange infinitesimal automorphism, even if it is not locally leafwise affine.

We finish this section by considering a more general structure.
\begin{defin} \label{def22} {\rm Let $(M,S)$ be a tangent
manifold. A {\it locally conformal Lagrange structure} on $M$ is a
maximal open covering $M=\cup U_\alpha$ with local regular
Lagrangians $ \mathcal{L}_\alpha$ such that, over the
intersections $U_\alpha\cap U_\beta$, the local Lagrangian metrics
satisfy a relation of the form
\begin{equation}\label{conformality}
g_{\mathcal{L}_\beta}=f_{(\alpha\beta)} g_{\mathcal{L}_\alpha},
\end{equation} where $ f_{(\alpha\beta)}>0$ are foliated functions. A tangent
manifold endowed with this type of structure is a {\it locally
conformal Lagrange space or manifold}.}\end{defin}

Clearly, condition (\ref{conformality}) is equivalent with the
transition relations \begin{equation}\label{conform2}
\mathcal{L}_\beta=f_{(\alpha\beta)}
\mathcal{L}_\alpha+a(\varphi_{(\alpha\beta)})+b_{(\alpha\beta)},\end{equation}
where the last two terms are like in (\ref{Lrel}). On the other
hand, $ \{\ln{ f_{(\alpha\beta)}}\}$ is a $\Phi$-valued
$1$-cocycle, and may be written as $\ln
f_{(\alpha\beta)}=\psi_\beta-\psi_\alpha$ where $\psi_\alpha$ is a
differentiable function on $U_\alpha$ (which may be assumed
projectable only if the cocycle is a coboundary). Accordingly the
formula \begin{equation}\label{gconformal}g|_{U_\alpha}=
e^{-\psi_\alpha}g_{ \mathcal{L}_\alpha}\end{equation} defines a
global transversal metric of the vertical foliation, which is
locally conformal with local Lagrange metrics. As a matter of
fact, we have
\begin{prop}\label{propos26} Let $(M^{2n},S)$ be a tangent
manifold, and $n>1$. Then, $M$ is locally conformal Lagrange iff
$M$ has a global transversal metric $g$ of the vertical foliation,
which is locally conformal with local Lagrange metrics.
\end{prop} \noindent {\bf Proof.} We still have to prove that
the existence of the metric $g$ that satisfies
(\ref{gconformal}) implies (\ref{conformality}), which is clear,
except for the fact that the functions $
f_{(\alpha\beta)}=e^{\psi_\beta-\psi_\alpha}$ are projectable.
This follows from the Lagrangian character of the metrics $g_{
\mathcal{L}_\alpha}$. Indeed, with the usual local coordinates
$(x^i,y^i)$, the symmetry of the derivative tensors $C$ of $g_{
\mathcal{L}_\alpha},g_{ \mathcal{L}_\beta}$ implies
$$\frac{\partial f_{(\alpha\beta)}}{\partial y^k}(g_{
\mathcal{L}_\alpha})_{ij}= \frac{\partial
f_{(\alpha\beta)}}{\partial y^i}(g_{ \mathcal{L}_\alpha})_{kj},$$
and a contraction by $(g_{ \mathcal{L}_\alpha})^{ij}$ yields
$\partial f_{(\alpha\beta)}/\partial y^k=0$. Q.e.d.

The cohomology class $\eta=[\ln{f_{(\alpha\beta)}}]\in
H^1(M,\Phi)$ will be called the {\it complementary class} of the
metric $g$, and the locally conformal Lagrange metric $g$ is a
locally Lagrange metric iff $\eta=0$. Indeed, if $\eta=0$, we may
assume that the functions $\psi_\alpha$ are foliated, and the
derivative tensor $C$ of $g=e^{-\psi_\alpha} g_{
\mathcal{L}_\alpha}$ is completely symmetric.

Furthermore, using a normal bundle $N\mathcal{V}$ and the leafwise
version of the de Rham theorem, the complementary class may be
seen as the $d''$-cohomology class of the global, $d''$-closed
{\it complementary form} $\tau$ obtained by gluing up the local
forms $\{d''\psi_\alpha\}$. In particular, Lemma
\ref{contractibility} and Proposition \ref{propos23} imply that
any locally conformal Lagrange metric $g$ of a tangent bundle must
be a locally, therefore, a globally Lagrange metric.
\begin{example}\label{example25} {\rm Consider the Hopf manifold
$H^{2n}$ of Example \ref{example4}. The local functions $
\sum_{i=1}^n(y^i)^2$ define a locally conformal Lagrange structure
on $H^{2n}$, and
$$g=\frac{\sum_{i=1}^n(y^i)^2}{\sum_{i=1}^n[(x^i)^2+(y^i)^2]}$$
is a corresponding global metric, which, with the previously used
notation, corresponds to $$\psi_\alpha=\ln{\{\sum_{i=1}^n[(x^i)^2+
(y^i)^2]\}}.$$  The corresponding complementary form is
$$\tau=\frac{2\sum_{i=1}^n y^idy^i}{\sum_{i=1}^n[(x^i)^2+
(y^i)^2]}.$$} \end{example}
\begin{prop}\label{propos27} Let $(M,S)$ be a tangent manifold and
$g$ a global transversal metric of the vertical foliation $
\mathcal{V}$ of $S$. Then, $g$ is locally conformal Lagrange iff
there exists a $d''$-closed form $\tau$ of type $(0,1)$ such that
the tensor $\tilde C=C-(\tau\circ S)\otimes g$, where $C$ is the
derivative tensor of $g$, is a completely symmetric
tensor.\end{prop} \noindent {\bf Proof.} Define $\tilde
g=e^{-\psi_\alpha}g$, where $\tau|_{U_\alpha}=d''\psi_\alpha$ for
a covering $M=\cup U_\alpha$. Then, $e^{-\psi_\alpha}\tilde C$ is
the derivative tensor of $\tilde g$, and the result follows from
Proposition \ref{propos25}. Q.e.d.
\section{Transversal Riemannian geometry}
The aim of this section is to give an index free presentation of the
connections used in Finsler and Lagrange geometry \cite{{Bao},{M1},{M2}},
while also extending these connections to tangent manifolds.

Let $(M,S)$ be a tangent manifold and $g$ a metric of the
transversal bundle of the vertical foliation $ \mathcal{V}$ ($
T\mathcal{V}=im\,S$). (The metrics which we consider are non
degenerate, but may be indefinite.) We do not get many interesting
differential-geometric objects on $M$, unless we fix a normal
bundle $N\mathcal{V}$, also called the {\it horizontal bundle},
i.e., we decompose
\begin{equation} \label{normal} TM=N\mathcal{V}\oplus
T\mathcal{V}.\end{equation} We will say that $N\mathcal{V}$ is a
{\it normalization}, and $(M,S,N\mathcal{V})$ is a {\it normalized
tangent manifold}. Where necessary, we shall use the local bases
(\ref{bases}), (\ref{cobases}). The projections on the two terms
of (\ref{normal}) will be denoted by $p_N$, $p_T$, respectively,
and $P=p_N-p_T$ is an almost product structure tensor that has the
horizontal and vertical distribution as $\pm
1$-eigendistributions, respectively.

For a normalized tangent manifold, the following facts are well
known: {\it i)} $S|_{N\mathcal{V}}$ is an isomorphism
$Q:N\mathcal{V}\rightarrow T\mathcal{V}$, {\it ii)} $S=Q\oplus0$,
{\it iii)} $S'=0\oplus Q^{-1}$ is an almost tangent structure,
{\it iv)} $F=S'+S$ is an almost product structure, {\it v)}
$J=S'-S$ is an almost complex structure on $M$.

On a normalized tangent manifold $(M,S,N\mathcal{V})$, a
pseudo-Riemannian metric $\gamma$ is said to be a {\it compatible
metric} if the subbundles $T\mathcal{V},N\mathcal{V}$ are
orthogonal with respect to $\gamma$ and \begin{equation}
\label{compatiblemetric} \gamma(SX,SY)=\gamma(X,Y),\hspace{5mm}
\forall X,Y\in\Gamma N\mathcal{V}.\end{equation} It is easy to see
that these conditions imply the compatibility of $\gamma$ with the
structures $J$ and $F$ i.e.,
\begin{equation}
\label{compatiblemetric2} \gamma(JX,JY)=\gamma(X,Y),\,
\gamma(FX,FY)=\gamma(X,Y), \hspace{5mm} \forall X,Y\in\Gamma
TM.\end{equation}

Furthermore, if $(M,S)$ is a tangent manifold and $\gamma$ is a
pseudo-Riemannian metric on $M$, we will say that $\gamma$ is
compatible with the tangent structure $S$ if the
$\gamma$-orthogonal bundle $N\mathcal{V}$ of $im\,S$ is a
normalization, and $\gamma$ is compatible for the normalized
tangent manifold $(M,S,N\mathcal{V})$.

The following result is obvious
\begin{prop} \label{proposit31} On a normalized, tangent
manifold, any transversal metric $g$ of the vertical foliation
defines a unique compatible metric $\gamma$, such that
$\gamma|_{N\mathcal{V}}=g$. \end{prop}

In what follows, we will refer at the metric $\gamma$ as the {\it
canonical extension} of the transversal metric $g$. On the other
hand, a pseudo-Riemannian metric $\gamma$ of a tangent manifold
$(M,S)$ which is the canonical extension of a locally Lagrange
metric $g$ will be called a {\it locally Lagrange-Riemann metric}.
This means that the restriction of $\gamma$ to the
$\gamma$-orthogonal subbundle $N\mathcal{V}$ of the vertical
foliation $ \mathcal{V}$ of $S$ is a locally Lagrange metric $g=
g_{\mathcal{L}_\alpha}$, and $\gamma$ is compatible with
$(M,S,N\mathcal{V})$ . Then, $(M,S,\gamma)$ will be called a {\it
locally Lagrange-Riemann manifold}. Notice that, since the induced
metric of $N\mathcal{V}$ is non degenerate, $N\mathcal{V}$ is a
normalization of the vertical foliation, and the compatibility
condition of the definition makes sense. Thus, any normalized
locally Lagrange space with the canonical extension $\gamma$ of
the Lagrange metric $g$ is a locally Lagrange-Riemann manifold,
and conversely.
\begin{example} \label{example32} {\rm The Euclidean metric
$\sum_{i=1}^n[(dx^i)^2+(dy^i)^2]$ is the canonical extension of
the locally Lagrange metric defined in Example \ref{example21} on the torus $
\mathbb{T}^{2n}$.}\end{example}
\begin{example}\label{example33} {\rm The metric
$$\sum_{i=1}^n(dx^i)^2+(dy)^2+\sum_{i=1}^n(dz^i-x^idy)^2+(dt)^2$$
is the canonical extension of the locally Lagrange metric defined
in Example \ref{example22} on $M(1,p)\times(
\mathbb{R}/\mathbb{Z})$.}\end{example}

Now, let $(M,S,N\mathcal{V},g)$ be a normalized tangent manifold
with a transversal metric of the vertical foliation $ \mathcal{V}$
and let $\nabla$ be the Levi-Civita connection of the canonical
extension $\gamma$ of $g$.

We are going to define a general connection that includes the
connections used in Finsler and Lagrange geometry
\cite{{Bao},{M1},{M2}}) as particular cases determined by specific
normalizations. This will be the so-called {\it second canonical
connection} $D$ of a foliated, pseudo-Riemannian manifold
$(M,\gamma)$, defined the following conditions \cite{V3}: i)
$N\mathcal{V}$ and $T\mathcal{V}$ are parallel, ii) the
restrictions of the metric to $N\mathcal{V}$ and $T\mathcal{V}$
are preserved by parallel translations along curves that are
tangent to $N\mathcal{V},T\mathcal{V}$, respectively, iii) the $
\mathcal{V}$-normal, respectively $ \mathcal{V}$-tangent,
component of the torsion $T_D(X,Y)$ vanishes if one of the
arguments is normal, respectively tangent, to $ \mathcal{V}$. This
connection is given by
\begin{equation} \label{secondc} \begin{array}{ll} D_{Z_1}Z_2
=p_N\nabla_{Z_1}Z_{2},&D_{Y_1}Y_2 =
p_T\nabla_{Y_1}Y_{2},\vspace{1mm}\\
D_{Y_1}Z_2 =p_N[{Y_1},Z_{2}],&D_{Z_1}Y_2 =p_T[{Z_1},Y_{2}],
\end{array} \end{equation} where $Y_1,Y_2\in\Gamma T\mathcal{V}$
and $Z_1,Z_2\in\Gamma N\mathcal{V}$. We will say that $D$ is the
{\it canonical connection}, and the connection induced by $D$ in
the normal bundle $N\mathcal{V}$, or, equivalently, in the
transversal bundle $\nu\mathcal{V}=TM/T\mathcal{V}$, will be
called the {\it canonical transversal connection}. The canonical,
transversal connection is a Bott (basic) connection \cite{Mol}.
The total torsion of the connection $D$ is not zero, namely one
has
\begin{equation} \label{torsionD} T_D(X,Y)=-p_T[p_NX,p_NY],
\hspace{5mm}\forall X,Y\in\Gamma TM. \end{equation}

\begin{prop}\label{propos33} Let $(M,S,g)$ be a locally
Lagrange manifold, and $\gamma$ the canonical extension of $g$.
Then, the derivative tensor field of $g$ has the following
expressions
\begin{equation} \label{CRiemann}
\begin{array}{l}C(X,Y,Z)=(D_{SX}g)(Y,Z)
=(D_{SX}\gamma)(Y,Z)\vspace{2mm}\\
=\gamma(\nabla_Y(SX),Z) +\gamma(Y,\nabla_Z(SX)),\end{array}
\end{equation} where $X,Y,Z\in\Gamma N\mathcal{V}$.
\end{prop} \noindent{\bf Proof.} Of course, in (\ref{CRiemann}),
$g$ is seen as a $2$-covariant tensor field on $M$ (see Section
2). First, we refer to the first two equalities (\ref{CRiemann}).
These are pointwise relations, hence, it will be enough to prove
these equalities for foliated cross sections of the normal bundle
$N\mathcal{V}$. Indeed, a tangent vector at a point can always be
extended to a projectable vector field on a neighborhood of that
point. But, in this case, the first and second equalities  are
straightforward consequences of the definitions of the tensor
field $C$ and of the connection $D$. Then, since $\nabla$ has no
torsion, (\ref{secondc}) implies
\begin{equation} \label{D-nabla}
D_{SX}Y=\nabla_{SX}Y-p_T\nabla_{SX}Y-p_N\nabla_Y(SX),
\end{equation}
and, also using $\nabla\gamma=0$, we get the required result.
Q.e.d.

The first two expressions of $C$ actually hold for any vector
fields $X,Y,Z\in\Gamma TM$.
\begin{corol} \label{corolar32} The canonical extension $\gamma$
of a transversal metric $g$ is a locally Lagrange-Riemann metric
iff one of the following two equivalent relations holds
\begin{equation}\label{simetria1}\begin{array}{rcl}
(D_{SX}\gamma)(Y,Z)&
=&(D_{SY}\gamma)(X,Z),\vspace{1mm}\\
\gamma(\nabla_Y(SX),Z) &+& \gamma(Y,\nabla_Z(SX))\vspace{1mm}\\ =
\gamma(\nabla_X(SY),Z)&+& \gamma(X,\nabla_Z(SY)),
\end{array}
\end{equation} where
$X,Y,Z\in\Gamma N\mathcal{V}$. \end{corol}
\begin{corol}\label{corolar31} On a tangent manifold, if $\gamma$
is a compatible pseudo-Rieman\-nian metric such that $\nabla S=0$,
then $\gamma$ is a projectable, locally Lagrange-Riemann metric.\end{corol}
\noindent {\bf Proof.} If $\nabla S=0$, the third equality
(\ref{CRiemann}) yields $C=0$, which is the characterization of this
type of metrics. Q.e.d.

Now we consider the curvature of $D$. The curvature is a tensor,
and it suffices to evaluate it pointwisely. For this reason,
whenever we need an evaluation of the curvature (as well as of any
other tensor) that involves vector fields, it will suffice to make
that evaluation on $\mathcal{V}$-projectable vector fields.
\begin{prop} \label{valoricurb} The curvature $R_D$
of the canonical connection has
the following properties
\begin{equation} \label{Bott1}
R_{D}(SX,SY)Z=0,\end{equation}
\begin{equation}\label{Bott2}R_{D}(SX,Y)Z=p_N[SX,D_YZ],\end{equation}
\begin{equation}\label{Bott3}
R_{D}(X,Y)(SZ)=-D_{SZ}(p_T[X,Y]),
\end{equation}
\begin{equation}\label{Bott4} R_{D}(SX,Y)Z=R_{D}(SX,Z)Y,
\end{equation}
for any foliated vector fields $X,Y,Z\in\Gamma N\mathcal{V}$.
Moreover, formulas (\ref{Bott1}), (\ref{Bott3}), and (\ref{Bott4})
hold for any arguments $X,Y,Z\in\Gamma N\mathcal{V}$.\end{prop}
\noindent{\bf Proof.} Equality (\ref{Bott1}) is in agreement with
the fact that $D$ is a Bott connection \cite{Mol}. Formulas
(\ref{Bott1})-(\ref{Bott3}) follow from (\ref{secondc}) and
(\ref{torsionD}). Formula (\ref{Bott4}) is a consequence of
(\ref{Bott2}). In the computation, one will take into account the
fact that for any foliated vector field $X\in\Gamma TM$ and any
vector field $Y\in\Gamma T\mathcal{V}$ one has $[X,Y]\in\Gamma
T\mathcal{V}$ \cite{Mol}. Q.e.d.
\begin{prop} \label{proposit30} For the canonical connection $D$,
the first Bianchi identity is equivalent to the following
equalities, where $X,Y,Z\in\Gamma N\mathcal{V}$
\begin{equation}\label{Bianchi0} \sum_{Cycl(X,Y,Z)}
R_D(SX,SY)(SZ)=0,
\end{equation}
\begin{equation}\label{Bianchi1} R_{D}(SX,Z)SY=R_{D}(SY,Z)SX,
\end{equation} \begin{equation}\label{Bianchi3} \sum_{Cycl(X,Y,Z)}
R_D(X,Y)Z=0.\end{equation} \end{prop} \noindent{\bf Proof.} Write
down the general expression of the Bianchi identity of a linear
connection with torsion (e.g., \cite{KN}) for arguments tangent
and normal to $ \mathcal{V}$. Then, compute using (\ref{secondc}),
(\ref{torsionD}) and projectable vector fields as arguments. The
fourth relation included in the Bianchi identity reduces to
(\ref{Bott3}). Q.e.d.
\begin{prop} \label{proposit39} For the canonical connection $D$,
the second Bianchi identity is equivalent to the following
equalities, where $X,Y,Z\in\Gamma N\mathcal{V}$,
\begin{equation} \label{Bianchi23}
\sum_{Cycl(X,Y,Z)}(D_{SX}R_{D})(SY,SZ)=0.\end{equation}
\begin{equation}\label{Bianchi21} (D_{SX}R_{D})(SY,Z)-(D_{SY}R_{D})(SX,Z)
=(D_ZR_D)(SX,SY),
\end{equation} \begin{equation} \label{Bianchi22}
(D_{X}R_{D})(Y,SZ)-(D_{Y}R_{D})(X,SZ)+(D_{SZ}R_{D})(X,Y)
\end{equation}$$=R_D(p_T[X,Y],SZ),$$
\begin{equation}\label{Bianchi20} \sum_{Cycl(X,Y,Z)}(D_{X}R_{D})(Y,Z)=
\sum_{Cycl(X,Y,Z)}R_D(p_T[X,Y],Z),\end{equation}
\end{prop}
\noindent{\bf Proof.} This is just a rewriting of the classical
second Bianchi identity \cite{KN} that uses (\ref{torsionD}).
Q.e.d.

Like in Riemannian geometry, we also define a covariant
curvature tensor \begin{equation} \label{covarcurb}
R_D(X,Y,Z,U)=\gamma(R_D(Z,U)Y,X),\hspace{5mm}X,Y,Z,U\in\Gamma TM.
\end{equation}
In particular, we have \begin{prop} \label{curb+C}
\begin{equation} \label{cocurb}
R_{D}(U,Z,SX,Y)=g([SX,D_YZ],U)\end{equation}
$$=(SX)(g(D_YZ,U))-C(X,D_YZ,U),$$ where the arguments
are foliated vector fields in $\Gamma N\mathcal{V}$,
and $g$ is seen as a tensor on $M$. \end{prop}

Formula (\ref{Bianchi3}) yields the Bianchi identity
\begin{equation}\label{Bianchicovar1}
\sum_{Cycl(X,Y,Z)}R_D(U,X,Y,Z)=0,\hspace{5mm}\forall
X,Y,Z\in\Gamma N\mathcal{V}. \end{equation} But, the other
Riemannian symmetries may not hold. Indeed, we have
\begin{prop} \label{antisym1-2} For any arguments $X,Y,Z,U\in
\Gamma N\mathcal{V}$ one has \begin{equation}\label{arg1-2}
R_D(X,Y,Z,U)+R_D(Y,X,Z,U)=(D_{p_T[Z,U]}\gamma)(X,Y)\end{equation}
$$=C(S'p_T[Z,U],X,Y).$$
\end{prop} \noindent{\bf Proof.} Express the equality
$$(ZU-UZ-[Z,U])(\gamma(X,Y))=0$$ for normal foliated arguments,
and use the transversal metric character of the canonical
connection $D$ and Proposition \ref{propos33}. Q.e.d.

\begin{prop}\label{symperechi} For any arguments $X,Y,Z,U\in
\Gamma N\mathcal{V}$ one has \begin{equation}\label{arg1-2,3-4}
R_D(X,Y,Z,U)-R_D(Z,U,X,Y)\end{equation}
$$=\frac{1}{2}\{C(S'p_T[Z,U],X,Y) - C(S'p_T[X,Y],Z,U)\}.$$
\end{prop} \noindent{\bf Proof.} Same proof as for
Proposition 1.1 of \cite{KN}, Chapter V. Q.e.d.

The other first and second
Bianchi identities may also be expressed in a covariant form.
From (\ref{covarcurb}) we get \begin{equation} \label{DRcovar}
(D_FR_D)(A,B,C,E)=\gamma((D_FR_D)(C,E)B,A)\end{equation}
$$+(D_F\gamma)(R_D(C,E)B,A),$$ where $(A,B,C,E,F\in \Gamma TM)$.
Accordingly, (\ref{Bianchi22}) yields
\begin{equation}\label{Bianchi21covar} (D_{SZ}R_D)(V,U,X,Y)
+(D_XR_D)(V,U,Y,SZ)-(D_YR_D)(V,U,X,SZ) \end{equation}
$$=(D_{SZ}\gamma)(R_D(X,Y)U,V) -(D_X\gamma)(p_N[SZ,D_YU],V)
+(D_Y\gamma)(p_N[SZ,D_XU],V), $$ (\ref{Bianchi20}) yields
\begin{equation}\label{Bianchi22covar}
\sum_{Cycl(X,Y,Z)}(D_XR_D)(V,U,Y,Z)=
\sum_{Cycl(X,Y,Z)}R_D(V,U,p_T[X,Y],Z)\end{equation}
$$-\sum_{Cycl(X,Y,Z)}(D_X\gamma)(R_D(Y,Z)U,V), $$ etc., where
$X,Y,Z,U,V\in \Gamma N\mathcal{V}$.
\begin{example}\label{exemplu1} {\rm On the torus $
\mathbb{T}^{2n}$ with the metric of Example \ref{example32}, the
usual flat connection is both the Levi-Civita connection and the
canonical connection $D$, and it has zero curvature. On the
manifold $M(1,p)\times( \mathbb{R}/\mathbb{Z})$ with the metric of
Example \ref{example33}, the connection that parallelizes the
orthonormal basis shown by the expression of the metric is not the
Levi-Civita connection, since it has torsion, but, it follows
easily that it has the characteristic properties of the canonical
connection $D$. Accordingly, we are in the case of a locally Lagrange-Riemann
manifold with a vanishing curvature $R_D$ and a non vanishing torsion $T_D$.}
\end{example}

\begin{prop} \label{proposit38} The Ricci curvature tensor
$\rho_D$ of the connection $D$ is given by the equalities
\begin{equation}\label{Ricci1} \rho_D(SX,SY)=
\sum_{i=1}^n<\vartheta^i,R_D(\frac{\partial}{\partial y^i},SX)SY>,
\end{equation} \begin{equation}
\label{Ricci2}\rho_D(SX,Y)= \sum_{i=1}^n<dx^i,p_N[D_{X_i}Y,SX]>,
\end{equation}
\begin{equation}\label{Ricci3} \begin{array}{rcl}
\rho_D(X,Y)&=&tr[Z\mapsto R_D(Z,X)Y]\vspace{1mm}\\
&=&\sum_{i=1}^n<dx^i,R_D(X_i,X)Y>,\end{array}\end{equation} where
$X,Y,Z\in\Gamma N\mathcal{V}$, and in (\ref{Ricci2}) $Y$ is projectable.
\end{prop} \noindent{\bf Proof.} The definition of the Ricci
tensor of a linear connection (e.g., \cite{KN}),
and the use of the bases (\ref{bases}) and (\ref{cobases}) yield
\begin{equation}\label{Ricci} \rho_D(X,Y)= \sum_{i=1}^n<dx^i,R_D(X_i,X)Y>
+ \sum_{i=1}^n<\theta^i,R_D(\frac{\partial}{\partial y^i},X)Y>.
\end{equation} Then, the results follow from (\ref{secondc})
and (\ref{Bott2}). Q.e.d.
\begin{rem} \label{curbscalara}
{\rm In view of (\ref{Ricci3}), we may speak of
$\kappa_D=tr\,\rho_D$ on $N\mathcal{V}$, and call it the {\it transversal
scalar curvature}.}\end{rem}

In the case of a normalized, bundle-type, tangent manifold
$(M,S,E,N\mathcal{V})$, with a compatible metric $\gamma$ ($E$ is
the Euler vector field), the curvature has some more interesting
features, which were studied previously in Finsler geometry
\cite{Bao}. These features follow from \begin{lemma}
\label{lemaDSZS'E} For any $Z\in\Gamma N\mathcal{V}$ one has
\begin{equation}\label{eqlemei} D_{SZ}(S'E)=Z.
\end{equation}\end{lemma} \noindent{\bf Proof.} $S'$ is the tensor defined
at the beginning of this section, and with local bundle-type
coordinates $(x^i,y^i)_{i=1}^n$ and bases (\ref{bases}), we have
$$SZ=\xi^i(x^j,y^j)\frac{\partial}{\partial y^i},
\;S'E=y^iX_i.$$ Now, (\ref{eqlemei}) follows
from (\ref{secondc}). Q.e.d.

Using Lemma \ref{lemaDSZS'E} one can prove
\begin{prop} \label{valoricurbura}
The curvature operator $R_D(X,Y)|_{N\mathcal{V}}$
$(X,Y\in\Gamma N\mathcal{V})$ is determined by its action
on $S'E$ and by  $R_D(V,SW)|_{N\mathcal{V}}$ where $V,W\in\Gamma N\mathcal{V}$.
\end{prop}
\noindent{\bf Proof.} Denote \begin{equation} \label{defr}
r(X,Y)=R_D(X,Y)S'E. \end{equation} The covariant derivative of
this tensor contains a term, which, in view of (\ref{eqlemei}), is
equal to $R_D(X,Y)Z$, and we get
\begin{equation} \label{ecuatiar} R_D(X,Y)Z=D_{SZ}(r(X,Y))
-r(D_{SZ}X,Y) \end{equation}
$$-r(X,D_{SZ}Y)-(D_{SZ}R_D)(X,Y)S'E.$$
Now, if the last term of (\ref{ecuatiar}) is expressed by means of
the Bianchi identity (\ref{Bianchi22}) one gets an expression of
$R_D(X,Y)Z$ in terms of $r$ and $R_D(V,SW)|_{N\mathcal{V}}$ for
various arguments $V,W$. Q.e.d.

Notice that, by (\ref{Bott2}), the computation of $R_D(V,SW)|_{N\mathcal{V}}$ on
normal arguments requires only a first order covariant derivative.

From Proposition \ref{valoricurbura} we also see that the
curvature values $R_D(U,Z,X,Y)$ $(X,Y,Z,U\in\Gamma N\mathcal{V})$
are determined by the values $R_D(U,S'E,X,Y)$ and of
$R_D(U,V,W,SK)$ for convenient normal arguments. Therefore, it
should be interesting to study manifolds where $R_D(U,S'E,X,Y)$
has a simple expression. If we fix a direction $span\{U\}$ and a
$2$-dimensional plan $\sigma=span\{X,Y\}$ $(U,X,Y\in\Gamma
N\mathcal{V})$, the formula \begin{equation} \label{Ucurbsect}
k_U(\sigma)=\frac{R_D(U,S'E,X,Y)}{\gamma(S'E,X)\gamma(U,Y) -
\gamma(S'E,Y)\gamma(U,X)}
\end{equation} defines an invariant, which we will call the $U$-{\it sectional
curvature} of $\sigma$. $k_U(\sigma)$ is independent of $U$ iff
\begin{equation} \label{Uindep} R_D(X,Y)S'E=k(\sigma)
[\gamma(S'E,X)Y-\gamma(S'E,Y)X], \end{equation}
where $k(\sigma)$ is a function of the point of $M$ and the plan $\sigma$ only.
Furthermore, if $k(\sigma)=f(x)$, $x\in M$, i.e. $k(\sigma)$ is pointwise constant,
(\ref{Uindep}) is a
natural simple expression of the transversal curvature tensor.

On the other hand, we can generalize the notion of {\it flag
curvature}, which is an important invariant in Finsler geometry
\cite{Bao}. Namely, a {\it flag} $\phi$ at a point $x\in M$ is a
$2$-dimensional plane $\phi\subseteq T_xM$ which contains the
vector $E_x$. Such a flag is $\phi= span\{E_x,X_x\}$, where
$X_x\in N_x\mathcal{V}$ is defined up to a scalar factor, and
following \cite{Bao}, the flag curvature is defined by
\begin{equation} \label{flagcurv}
k(\phi)=k(X)=\frac{R_D(X,S'E,X,S'E)}{g(S'E,S'E)g(X,X)-g^2(S'E,X)}.
\end{equation}

If $g$ is not positive definite, the flag curvature may take
infinite values.
\begin{prop} \label{propflag} The flag curvature $k$ is pointwise
constant iff \begin{equation} \label{constantflagc}
R_D(X,S'E,Y,S'E)=f[g(S'E,S'E)g(X,Y)-g(S'E,X)g(S'E,Y)],\end{equation}
where $f\in C^\infty(M)$. If
the $U$-sectional curvature is independent of $U$ and poinwise constant,
the flag curvature is pointwise constant too.
\end{prop} \noindent{\bf Proof.} For the first assertion, use
$k(X+Y)=k(X)=k(Y)$. The second follows because, if
$k(\sigma)=f(x)$, (\ref{Uindep}) implies \begin{equation}
\label{Rtriplu}
R_D(U,S'E,X,Y)=f(x)[\gamma(S'E,X)\gamma(Y,U)-\gamma(S'E,Y\gamma(X,U)],
\end{equation} which reduces to
(\ref{constantflagc}) for $Y=S'E$. Q.e.d.
\begin{rem} \label{obsFinsler} {\rm The curvature $R_D$ has more
interesting properties in the case of a bundle-type, locally Lagrange manifold
such that the metric tensor $g$ is homogeneous of degree zero
with respect to the coordinates $y^i$. The invariant characterization
of this situation is that
the derivative tensor $C$ is symmetric, and such that
\begin{equation} \label{eqobsFinsler} i(S'E)C=0.
\end{equation} Indeed, in this case, formulas (\ref{arg1-2}),
(\ref{arg1-2,3-4}), etc., yield simpler symmetry properties if one
of the arguments is $S'E$. The Finsler metrics satisfy the
homogeneity condition (\ref{eqobsFinsler}).} \end{rem}
\begin{rem}\label{alteconex} {\rm
On a locally Lagrange-Riemann manifold $(M,S,\gamma)$
there exist other geometrically interesting
connections as well. One such connection is
\begin{equation}\label{nablabar}
\nabla'_{X}Y=p_{N}(\nabla _X(p_{N}Y)) +
p_{T}(\nabla_{X}(p_{T}Y)).\end{equation} The connection $\nabla'$
preserves the vertical and horizontal distributions and the
metric, but has a non zero torsion. Then, we have the connections
$ ^C\hspace{-1mm}D,^C\nabla'$, which can be defined by using
formulas (\ref{secondc}), (\ref{nablabar}) with the Levi-Civita
connection $\nabla$ replaced by the {\it Chern connection}
$^C\nabla$ i.e., the $\gamma$-metric, $J$-preserving connection
that has a torsion with no component of $J$-type $(1,1)$
$(J=S'-S)$ \cite{KN}.}\end{rem}

We finish by recalling the well known fact \cite{Kern,M1,M2} that
global Finsler and Lagrange structures of tangent bundles have an
invariant normalization. This normalization may be defined as
follows.

Let $ \mathcal{L}$ be the global
Lagrangian function.Then the {\it energy function} \begin{equation}
\label{energy}\mathcal{E}_\mathcal{L}=E\mathcal{L}-\mathcal{L}
\end{equation} has a Hamiltonian vector field
$X_\mathcal{E}$ defined by
\begin{equation} \label{HamE} i(X_{{\mathcal
E}})\omega_\mathcal{L}=-d\mathcal{E}_\mathcal{L},
\end{equation} where $\omega_\mathcal{L}$ is the Lagrangian symplectic
form (\ref{formasimpl}), which turns out to be a second order
vector field. Accordingly, $L_{X_\mathcal{E}}S$ is an almost
product structure on $M$ (see Section 1), and
$N_\mathcal{E}\mathcal{V}= im\,H$, with $H$ defined by
(\ref{almprod}) is a canonical normal bundle of $ \mathcal{V}$.

A locally Lagrangian structure $\{\mathcal{L}_\alpha\}$ on a
bundle-type tangent manifold $(M,S,E)$ defines a global function
({\it second order energy})
\begin{equation} \label{secondenerg} \mathcal{E}'=E^2\mathcal{L}_\alpha
-E\mathcal{L}_\alpha, \end{equation} but, generally, it has no
global Hamiltonian vector field,
and, even if such a field exists, is may not be a second order vector field.\vspace{2mm}\\
\noindent{\it Acknowledgement}. Part of the work on this paper was done
during a visit at the Erwin Schr\"odinger International Institute for
Mathematical Physics, Vienna, Austria, October 1-10, 2002, in the
framework of the program ``Aspects of foliation theory in geometry, topology
and physics". The author thanks the organizers of the program, J. Glazebrook,
F. Kamber and
K. Richardson, the ESI and its director prof. P. Michor
for having made that visit possible.
\hspace*{7.5cm}{\small \begin{tabular}{l} Department of
Mathematics\\ University of Haifa, Israel\\ E-mail:
vaisman@math.haifa.ac.il \end{tabular}}
\end{document}